\documentclass[12pt]{amsart}
\usepackage{amssymb}
\usepackage{verbatim}
\usepackage{hyperref}
\usepackage{color}

\allowdisplaybreaks
\begin{document}

% define theorem environments
\newtheorem{theorem}{Theorem}    %[section]
\newtheorem{proposition}[theorem]{Proposition}
\newtheorem{conjecture}[theorem]{Conjecture}
\def\theconjecture{\unskip}
\newtheorem{corollary}[theorem]{Corollary}
\newtheorem{lemma}[theorem]{Lemma}
\newtheorem{sublemma}[theorem]{Sublemma}
\newtheorem{observation}[theorem]{Observation}
\theoremstyle{definition}
\newtheorem{definition}{Definition}
\newtheorem{notation}[definition]{Notation}
\newtheorem{remark}[definition]{Remark}
\newtheorem{question}[definition]{Question}
\newtheorem{questions}[definition]{Questions}
\newtheorem{example}[definition]{Example}
\newtheorem{problem}[definition]{Problem}
\newtheorem{exercise}[definition]{Exercise}

\numberwithin{theorem}{section} \numberwithin{definition}{section}
\numberwithin{equation}{section}

\def\earrow{{\mathbf e}}
\def\rarrow{{\mathbf r}}
\def\uarrow{{\mathbf u}}
\def\varrow{{\mathbf V}}
\def\tpar{T_{\rm par}}
\def\apar{A_{\rm par}}

\def\reals{{\mathbb R}}
\def\torus{{\mathbb T}}
\def\heis{{\mathbb H}}
\def\integers{{\mathbb Z}}
\def\naturals{{\mathbb N}}
\def\complex{{\mathbb C}\/}
\def\distance{\operatorname{distance}\,}
\def\support{\operatorname{support}\,}
\def\dist{\operatorname{dist}\,}
\def\Span{\operatorname{span}\,}
\def\degree{\operatorname{degree}\,}
\def\kernel{\operatorname{kernel}\,}
\def\dim{\operatorname{dim}\,}
\def\codim{\operatorname{codal}}
\def\trace{\operatorname{trace\,}}
\def\Span{\operatorname{span}\,}
\def\dimension{\operatorname{dimension}\,}
\def\codimension{\operatorname{codimension}\,}
\def\nullspace{\scriptk}
\def\kernel{\operatorname{Ker}}
\def\ZZ{ {\mathbb Z} }
\def\p{\partial}
\def\rp{{ ^{-1} }}
\def\Re{\operatorname{Re\,} }
\def\Im{\operatorname{Im\,} }
\def\ov{\overline}
\def\eps{\varepsilon}
\def\lt{L^2}
\def\diver{\operatorname{div}}
\def\curl{\operatorname{curl}}
\def\etta{\eta}
\newcommand{\norm}[1]{ \|  #1 \|}
\def\expect{\mathbb E}
\def\bull{$\bullet$\ }
\def\C{\mathbb{C}}
\def\R{\mathbb{R}}
\def\Rn{{\mathbb{R}^n}}
\def\Sn{{{S}^{n-1}}}
\def\M{\mathbb{M}}
\def\N{\mathbb{N}}
\def\Q{{\mathbb{Q}}}
\def\Z{\mathbb{Z}}
\def\F{\mathcal{F}}
\def\L{\mathcal{L}}
\def\S{\mathcal{S}}
\def\supp{\operatorname{supp}}
\def\dist{\operatorname{dist}}
\def\essi{\operatornamewithlimits{ess\,inf}}
\def\esss{\operatornamewithlimits{ess\,sup}}
\def\xone{x_1}
\def\xtwo{x_2}
\def\xq{x_2+x_1^2}
\newcommand{\abr}[1]{ \langle  #1 \rangle}
\newcommand{\f}{\frac}

\newcommand{\Norm}[1]{ \left\|  #1 \right\| }
\newcommand{\set}[1]{ \left\{ #1 \right\} }
\def\one{\mathbf 1}
\def\whole{\mathbf V}
\newcommand{\modulo}[2]{[#1]_{#2}}

\def\scriptf{{\mathcal F}}
\def\scriptg{{\mathcal G}}
\def\scriptm{{\mathcal M}}
\def\scriptb{{\mathcal B}}
\def\scriptc{{\mathcal C}}
\def\scriptt{{\mathcal T}}
\def\scripti{{\mathcal I}}
\def\scripte{{\mathcal E}}
\def\scriptv{{\mathcal V}}
\def\scriptw{{\mathcal W}}
\def\scriptu{{\mathcal U}}
\def\scriptS{{\mathcal S}}
\def\scripta{{\mathcal A}}
\def\scriptr{{\mathcal R}}
\def\scripto{{\mathcal O}}
\def\scripth{{\mathcal H}}
\def\scriptd{{\mathcal D}}
\def\scriptl{{\mathcal L}}
\def\scriptn{{\mathcal N}}
\def\scriptp{{\mathcal P}}
\def\scriptk{{\mathcal K}}
\def\frakv{{\mathfrak V}}

\begin{comment}
\def\scriptx{{\mathcal X}}
\def\scriptj{{\mathcal J}}
\def\scriptr{{\mathcal R}}
\def\scriptS{{\mathcal S}}
\def\scripta{{\mathcal A}}
\def\scriptk{{\mathcal K}}
\def\scriptp{{\mathcal P}}
\def\frakg{{\mathfrak g}}
\def\frakG{{\mathfrak G}}
\def\boldn{\mathbf N}
\end{comment}

\title[Mutlilinear   square functions]{Certain Multi(sub)linear square functions}

%title of paper and the running head option

\author[Loukas Grafakos]{Loukas Grafakos} %first author's name and the running head option

\author[Sha He]{Sha He}%second author's name and the running head option

\author[Qingying Xue]{Qingying Xue}%second author's name and the running head option

%%%%%%%%%%%%%%% footnote %%%%%%%%%%%%%%%%
%In case \subjclass[2000] command is not effective
%(or the version of amsart.cls is old), write as follows instead:
%\renewcommand{\thefootnote}{\fnsymbol{footnote}}
%\footnote[0]{2000\textit{ Mathematics Subject Classification}.
%Primary 00; Secondary 00.}
%
\thanks{The first author was supported partly by the Simons Foundation. The third author was supported partly by NSFC
(No. 11471041), the Fundamental Research Funds for the Central Universities (No. 2014KJJCA10) and NCET-13-0065.} % \\ \indent Corresponding author: }
%%%%%%%%%%%% Authors addresses %%%%%%%%%%%%%
\address{Loukas Grafakos\\Department of Mathematics\\University of Missouri\\Columbia, MO, 65211\\
USA
}
\email{grafakosl@missouri.edu}

\address{She He\\School of Mathematical
Sciences\\Beijing Normal University\\Laboratory of Mathematics and
Complex Systems\\Ministry of Education\\Beijing, 100875\\
P. R. China}
\email{amyhesha@163.com}

\address{Qingying Xue\\School of Mathematical
Sciences\\Beijing Normal University\\Labo-ratory of Mathematics and
Complex Systems\\Ministry of Education\\ Beijing \linebreak 100875\\
P. R. China
}
\email{qyxue@bnu.edu.cn}

\begin{abstract}
Let $d\ge 1, \ell\in\Z^d$, $m\in \mathbb Z^+$ and $\theta_i$, $i=1,\dots,m $ are fixed, distinct and nonzero real numbers.
We show that the   $m$-(sub)linear version below of the Ratnakumar and Shrivastava \cite{RS1} Littlewood-Paley square function
$$
T(f_1,\dots , f_m)(x)=\Big(\sum\limits_{\ell\in\Z^d}|\int_{\mathbb{R}^d}f_1(x-\theta_1 y)\cdots f_m(x-\theta_m y)e^{2\pi i \ell \cdot y}K (y)dy|^2\Big)^{1/2}
$$
 is bounded from  $L^{p_1}(\mathbb{R}^d) \times\cdots\times L^{p_m}(\mathbb{R}^d) $ to $L^p(\mathbb{R}^d) $ when $2\le p_i<\infty$ satisfy $1/p=1/p_1+\cdots+1/p_m$ and  $1\le p<\infty$. Our proof is based on a modification  of an inequality of
 Guliyev and Nazirova \cite{GN} concerning multilinear convolutions.
\end{abstract}

\maketitle
%------------------------------------------------------------------------
\section{Introduction and main results}
Motivated by the study  of the  bilinear Hilbert transform
\begin{equation}\label{BHT}
H(f,g)(x)=p.v.\frac{1}{\pi}\int_{\R}f(x-y)g(x+y)\frac{dy}{y},
\end{equation}
by Lacey and Thiele \cite{LT1}, \cite{LT2}, a substantial amount of work has been produced in the area of multilinear singular integral and multiplier operators. In this paper we study certain kinds of multi(sub)linear square functions.

%Grafakos and Li \cite{GL} and Li \cite{L}.

We introduce a     bilinear operator which is closely related to the bilinear Hilbert transform.
Let $w$ be a cube in $\R^d$, $d\geq 1$. Then for functions $f,g$ in the Schwartz class $\mathcal{S}(\R^d)$, we define the bilinear operator $S_w$ associated with the symbol $\chi_w(\xi-\eta)$ as follows:
$$
S_w(f,g)(x)=\int_{\R^d}\int_{\R^d}\hat{f}(\xi)\hat{g}(\eta)\chi_w(\xi-\eta)e^{2\pi i x\cdot(\xi+\eta)}d\xi d\eta.
$$
We notice that when $d=1$ and $w$ is the characteristic function of a half plane in $\mathbb R^2$, then
$S_w$ is related to the bilinear Hilbert transform in \eqref{BHT}.

Let $\{w_l\}_{l\in\Z^d}$ be a sequence of disjoint cubes in $\R^d$. Let $S_{w_l}$ be the bilinear operator associated with the symbol $\chi_{w_l}(\xi-\eta)$ as defined above. Then for $f,g\in\mathcal{S}(\R^d)$, the {\it non-smooth bilinear Littlewood-Paley square function} associated with the sequence $\{w_l\}_{l\in\Z^d}$ is defined as
\begin{equation}\label{BSSF}
S(f,g)(x)=\Big(\sum\limits_{l\in\Z^d}|S_{w_l}(f,g)(x)|^2\Big)^{1/2}.
\end{equation}
A   {\it smooth } version of the
bilinear Littlewood-Paley square function   in \eqref{BSSF} can be obtained if the
characteristic function of the cube $w$ is replaced by a smooth bump adapted to $w$. Precisely, we let
\begin{equation}\label{BSSF2}
T(f,g)(x)=\Big(\sum\limits_{l\in\Z^d}|T_{\phi_l}(f,g)(x)|^2\Big)^{1/2},
\end{equation}
where $T_{\phi_l}$ is a bilinear operator associated with the smooth function $\phi_l$ whose Fourier transform is supported in $ w_l$, i.e.,
$$
T_{\phi_l}(f,g)(x)=\int_{\R^d}\int_{\R^d}\widehat{f}(\xi)\widehat{g}(\eta)\widehat{\phi_l}(\xi-\eta)e^{2\pi ix\cdot(\xi+\eta)}d\xi d\eta.
$$

A particular case of the bilinear Littlewood-Paley square function $T$ in \eqref{BSSF2} was   studied by Lacey \cite{L} who proved the following: \medskip\\
\noindent\textbf{Theorem A} (\cite{L}).
Let $\phi$ be a smooth function on $\R^d$ such that $\hat{\phi}$ is supported in the unit cube of $\R^d$. For $l\in\Z^d$, let $\widehat{\phi_l}$ be the function defined by $\widehat{\phi_l}(\xi)=\widehat{\phi}(\xi-l)$. Then for $2\leq p,q\leq\infty$ with $1/p+1/q=1/2$, there is a constant $C$ such that for $f,g\in\mathcal{S}(\R^d)$,
$$
\|T(f,g)\|_{L^2(\R^d)}\leq C\|f\|_{L^p(\R^d)}\|g\|_{L^q(\R^d)}.
$$

\medskip
In addition, Lacey raised two questions related to this:
(a) Does  {Theorem A} hold  for $r\neq 2$?
(b) Are there analogues for  non-smooth square functions?
The boundedness of non-smooth square functions on $\R\times \R$   strengthens the boundedness of bilinear Hilbert transform; in this direction and in the spirit of question (b), Bernicot \cite{B}, Ratnakumar and Shrivastava \cite{RS2} have provided some interesting results. As far as  question (a),  answers were provided by   Bernicot and  Shrivastava \cite{BS}, Mohanty and Shrivastava \cite{MS}.
The proofs of these results were based on  rather complicated  time-frequency analysis. However, Ratnakumar and Shrivastava \cite{RS1}   provided  a proof for the boundedness of a smooth bilinear square functions, which is based on more elementary techniques. Motivated by the work of \cite{RS1} and the increasing interest in multilinear operators, the aim of this paper is to study the $L^p$ boundedness properties of smooth $m$-(sub)linear Littlewood-Paley square functions.

\begin{definition}
We define an {\it $m$-(sub)linear Littlewood-Paley square function} as:
\begin{equation}\label{kkdfjd}
T(f_1,\dots,f_m)(x)=\Big(\sum\limits_{\ell\in\Z^d}|T_{\ell}(f_1,\dots,f_m)(x)|^2\Big)^{1/2},
\end{equation}
where
\begin{equation}\label{eq1}
T_{\ell}(f_1,\dots,f_m)(x)=\int_{\R^d}f_1(x-\theta_1 y)\cdots f_m(x-\theta_m y)K_{\ell}(y)dy{ .}
\end{equation}
Here $\theta_j$, $j=1,\dots,m$ are fixed, distinct and nonzero real numbers, and $K_{\ell}$ are integrable  functions that satisfy  certain additional size conditions.
\end{definition}

Expressing each $T_\ell$ in $m$-linear Fourier multiplier form, we write
\begin{align*}
&T_{\ell}(f_1,\dots,f_m)(x)=\int_{(\R^d)^m}\Big[\prod_{i=1}^m( \widehat{f_i}(\xi_i)
e^{2\pi i x\cdot\xi_i})\Big]\widehat{K_{\ell}}(\theta_1\xi_1+\cdots+\theta_m\xi_m)d\xi_1\cdots d\xi_m.
\end{align*}
Then $T_\ell$ is associated with the multiplier (or symbol)
$\widehat{K_{\ell}}(\theta_1\xi_1+\cdots+\theta_m\xi_m)$.

We use the notation $r'=r/(r-1)$ for the dual exponent of $r\in [1,\infty]$ with $1'=\infty$ and $\infty'=1$. Our main result is the following:
\begin{theorem}\label{thm1}
Let $1\le p<\infty$ and let $K$ be a measurable function on $\R^d$ which satisfies
\begin{equation}\label{kkjrff}
B_p=\sum\limits_{u\in \Z^d}\big\|\chi_{Q_u}K\big\|_{L^{p' }(\R^d)}<\infty  \qquad \textup{when $1\le p<2$}
\end{equation}
 and
$$
B_2=\sum\limits_{u\in \Z^d}\big\|\chi_{Q_u}K\big\|_{L^2(\R^d)} <\infty  \qquad \textup{when $2\le p<\infty$,}
$$
where $Q_u=u+[0, 1)^d$ for $u \in \Z^d$.
%$\sum\limits_{u\in \Z^d}=\sum\limits_{u_1\in \Z}\cdots\sum\limits_{u_d\in \Z}$.
For $\ell\in \Z^d$   define $ {K_{\ell}}(x)= {K}(x )e^{2\pi i  x \cdot  \ell }$ and let $T_{\ell}$ be the $m$-linear Fourier multiplier     in \eqref{eq1}, where the $\theta_j$ are nonzero and distinct.  Then for $2\le p_j<\infty$, $j=1,\dots,m$,
%$(1<)\frac{2r}{3r-2}<p<\infty$,
%$2<q<\infty$, $1< r< 2$ satisfying
%\frac{1}{q}+1=\frac{1}{p}+\frac{1}{r}\, , \qquad
there exists a positive constant $C=C(d, \theta_j  ,p_j)$ such that  the square function $T$ in \eqref{kkdfjd} satisfies
\begin{equation}\label{1.2}
\|T(f_1,\dots,f_m)\|_{L^p(\R^d)}\leq BC\prod\limits_{j=1}^m\|f_j\|_{L^{p_j}(\R^d)}
\end{equation}
for all functions $f_j$ in $L^{p_j}(\R^d)$, where $p$ and $p_j$ are related via
\begin{equation}\label{RES}
 \frac{1}{p}=\sum\limits_{j=1}^m\frac{1}{p_j}
\end{equation}
and $B$ in \eqref{1.2} is either $B_2$ if $2\le p<\infty  $ or $B_p$ when $1\le p<2$.
\end{theorem}

The case $m=2$ of this theorem was obtained  by Ratnakumar and Shrivastava \cite{RS1}.
The new ingredient of this note is the extension of this result to the case where $m\ge 3$, where a multilinear version of Young's inequality is needed. Such an inequality was obtained by
Guliyev and Nazirova \cite{GN} for a certain range of exponents which is not sufficient for our work and for this reason, we
provide an extension of  this paper.   To state the  multilinear Young  inequality of Guliyev and Nazirova, we fix $\theta_i$ nonzero distinct real numbers and we define
$$
\vec{f}\otimes g(x)=\int_{\R^d}f_1(x-\theta_1y)\cdots f_m(x-\theta_my)g(y)dy{ ,}
$$
where $\vec{f}=(f_1,\dots , f_m)$ and $f_j$, $g$ are measurable functions so that the preceding integral converges.
We denote by $L^{p,\infty}(\R^d)$ the weak $L^p$ space of all measurable functions $f$ such that
$
\|f\|_{L^{p,\infty}}=\sup\limits_{\alpha>0}\big\{\alpha|\{x\in\R^d:|f(x)|>\alpha\}|^{1/p}\big\}<\infty.
$

\begin{theorem}\label{lem1} (\cite{GN})
Assume that $1<r <\infty$,  $1< p_j \le \infty$, $i=1,\dots,m$,   $1/p=1/{p_1}+\cdots+1/{p_m}$,  $1/q+1=1/p+1/r$, and
$\frac{r'}{1+r'}\le  p<r'$ (equivalently $1\le q<\infty$). Then for $g\in L^{r,\infty}(\R^d)$ and
$f_j\in L^{p_j}(\R^d) $  we have $\vec{f}\otimes g\in L^q(\R^d)$ with    norm inequality
$$
\|\vec{f}\otimes g\|_{L^q(\R^d)}\leq C(d,\theta_j,p_j,r )\prod\limits_{i=1}^m\|f_i\|_{L^{p_i}(\R^d)}\|g\|_{L^{r,\infty}(\R^d)}.
$$
 \end{theorem}

Unfortunately,  the proof of Theorem  \ref{lem1} in \cite{GN} in the  range    $\frac{r'}{r'+1} \le p\le 1$
 (equivalently $1\le q\le r$) contains an inaccurate deduction (and the assertion is incorrect in this case).
In the proof of Theorem \ref{thm1} we need the case where $\frac{r'}{1+r'}=  p $ (equivalently $q=1$) and for this reason we
provide an analogue in which the weak norm of $g$ is replaced by a strong norm.
So we fix the inaccurate deduction in \cite{GN} via the following result:

 \begin{theorem}\label{lem11}
Assume that $1<r \le  \infty$,  $1/p=1/{p_1}+\cdots+1/{p_m}$, $1\le p_j \le \infty$, $j=1,\dots,m$,   $1/q+1=1/p+1/r$, and
$\frac{r'}{1+r'}\le  p\le r'$ (equivalently $1\le q \le \infty$). Then for some
constant $C(d,\theta_j,p_j,r )$  the following inequalities are valid:\\
(a) if $1<r<\infty$, then when $1<   p<r'$ (equivalently $r< q<\infty$) we have
$$
\|\vec{f}\otimes g\|_{L^q(\R^d)}\leq C(d,\theta_j,p_j,r )\prod\limits_{i=1}^m\|f_i\|_{L^{p_i}(\R^d)}\|g\|_{L^{r,\infty}(\R^d)}.
$$
(b) if  $1<r\le \infty$, then when $\frac{r'}{1+r'}\le  p\le 1$ (equivalently $1\le q\le r$) we have
$$
\|\vec{f}\otimes g\|_{L^q(\R^d)}\leq C(d, \theta_j,p_j ,r)\prod\limits_{i=1}^m\|f_i\|_{L^{p_i}(\R^d)}\|g\|_{L^{r }(\R^d)}.
$$
\end{theorem}

 \begin{remark}
 We notice that  in Theorem \ref{lem11} we are assuming that $p_j\ge 1$ for all $j$
whereas in \cite{GN} it is assumed that
$p_j>1$ for all $j$.   This small detail, i.e., the case $p_j=1$ for some $j$, does not appear in case (a), since  we
always have $p\le p_j$ for all $j$ and in case (a) we have $p >1$. Thus the case where some $p_j$ equals $1$
is relevant only  in case (b), which is treated in the next section.
  \end{remark}

 \section{The proof of Theorem \ref{lem11}}

The case where $\frac{r'}{r'+1} \le p\le 1$, i.e.,
case (a) of Theorem \ref{lem11} is contained in \cite{GN}.
So we focus attention in  case (b)   of Theorem \ref{lem11}, i.e., when
$\frac{r'}{1+r'}\le  p\le r'$ (equivalently $1\le q\le r$). Below we assume that $1<r\le \infty$.

First we consider the case $p=1$ when $q=r$. Since $1/p_1+\cdots+1/p_m=1$,
  making use of H\"older's inequality and Fubini's theorem, we write
\begin{align*}
&\|\vec{f}\otimes g\|_{L^r(\R^d)}\\
&\leq \Big(\int_{\R^d}\Big(\int_{\R^d}|f_1(x-\theta_1 y)\cdots f_m(x-\theta_m y)g(y)|dy\Big)^rdx\Big)^{\f 1 r}\\
&\leq \Big(\int_{\R^d}\Big\{\Big(\int_{\R^d}|g(y)|^r |f_1(x-\theta_1 y)\cdots f_m(x-\theta_m y)|dy\Big)^{\f 1 r}\\
&\ \ \ \qquad \qquad\qquad \qquad \qquad \Big(\int_{\R^d}|f_1(x-\theta_1 y)\cdots f_m(x-\theta_m y)|dy\Big)^{\f 1 {r'}}\Big\}^rdx\Big)^{\f 1 r}\\
&\leq \bigg(\prod_{i=1}^m |\theta_i|^{-\f{d}{p_i r'}}\bigg)\Big(\int_{\R^d} \int_{\R^d}|g(y)|^r |f_1(x-\theta_1 y)\cdots f_m(x-\theta_m y)|dy\,  dx\Big)^{\f 1 r}\prod\limits_{i=1}^m\|f_i\|_{L^{p_i}}^{\f 1{r'}}\\
&= \bigg(\prod_{i=1}^m |\theta_i|^{-\f{d}{p_i r'}}\bigg) \Big(\int_{\R^d}|g(y)|^r\int_{\R^d}|f_1(x-\theta_1 y)\cdots f_m(x-\theta_m y)|dx\, dy\Big)^{\f 1 r}\prod\limits_{i=1}^m\|f_i\|_{L^{p_i}}^{\f 1{r'}}\\
&\leq \bigg(\prod_{i=1}^m |\theta_i|^{-\f{d}{p_i r'}}\bigg)\|g\|_{L^r}\prod\limits_{i=1}^m\|f_i\|_{L^{p_i}}^{\f 1r}\prod\limits_{i=1}^m\|f_i\|_{L^{p_i}}^{\f 1{r'}}\\
&=\bigg(\prod_{i=1}^m |\theta_i|^{-\f{d}{p_i r'}}\bigg)\|g\|_{L^r}\prod\limits_{i=1}^m\|f_i\|_{L^{p_i}},
\end{align*}
where in the second inequality we apply H\"older's inequality with respect to the measure $\prod\limits_{i=1}^m|f_i(x-\theta_i y)|dy$ to the function $y\mapsto|g(y)|$ and 1 with exponents $r$ and $r'$,  respectively. Notice that this proof also
works in the case where $r=\infty$ with a small notational modification.

Next we consider the case $p=\frac{r'}{r'+1}$, which is equivalent to $q=1$. In this case $r$ could be equal to $\infty$ and thus  $r'$ could take the value   $1$ when $p=1/2$.
Let us give the proof first in the case where $m=3$, where the notation is simpler. In this case when $r'>1$
notice that $1/p_1+1/p_2+1/p_3=1+1/r'<2$, thus in the  space  with coordinates $(x_1,x_2,x_3)$, the intersection of the plane $x_1+x_2+x_3=1+1/r'$ and of the boundary of the unit cube $[0,1]^3$
is a closed hexagon with vertices
at six points $V_1^{(3)}=(1/r',0,1)$, $V_2^{(3)}=(1/r',1,0)$, $V_3^{(3)}=(1,0,1/r')$, $V_4^{(3)}=(0,1,1/r')$, $V_5^{(3)}=(1,1/r',0)$, $V_6^{(3)}=(0,1/r',1)$.
In the case where $r=\infty$, i.e., $r'=1$,   the intersection of boundary of the
unit cube $[0,1]^3$ with the plane $x_1+x_2+x_3= 2$ is a closed triangle with vertices
the points $W_1^{(3)}=(0,1,1)$, $W_2^{(3)}=(1,0,1)$ and $W_3^{(3)}=(1,1,0)$.

We assert that $\vec{f}\otimes g$ maps $L^{p_1} \times L^{p_2}\times L^{p_3}\times L^r$ into $L^1$ when $(1/p_1,1/p_2,1/p_3)$ lies in the closed convex hull of the  six points $V_j^{(3)}$ when $r<\infty$ (or the three points $W_1^{(3)}$, $W_2^{(3)}$, $W_3^{(3)}$ when $r=\infty$).   This assertion will be a consequence of the fact that boundedness holds at the   vertices via multilinear complex interpolation (in particular by applying the theorem in  Zygmund  \cite[Chapter XII, (3.3)]{Zyg} or in Berg and L\"ofstrom   \cite[Theorem 4.4.2]{BL} or   the main result in \cite{GrMa}).

 We only prove boundedness for point $V_1^{(3)}$, since the method is similar for the remaining vertices. We have
\begin{align*}
&\|\vec{f}\otimes g\|_{L^1(\R^d)}\\
&\leq \int_{\R^d}\int_{\R^d}|f_1(x-\theta_1 y)||f_2(x-\theta_2 y)||f_3(x-\theta_3 y)||g(y)|dydx\\
& =\int_{\R^d}|g(y)|\int_{\R^d}|f_1(x-\theta_1 y)||f_2(x-\theta_2 y)||f_3(x-\theta_3 y)|dxdy\\
&\leq \Big(\int_{\R^d}|g(y)|^rdy\Big)^{ \f 1r}\Big(\int_{\R^d}\Big(\int_{\R^d}|f_1(x-\theta_1 y)||f_2(x-\theta_2 y)||f_3(x-\theta_3 y)|dx\Big)^{r'}dy\Big)^{\f 1{r'}}\\
&\leq \|g\|_{L^r}\|f_2\|_{L^\infty}\Big(\int_{\R^d}\Big(\int_{\R^d}|f_1(t)|
|f_3(t-(\theta_3-\theta_2)y)|dt\Big)^{r'}dy\Big)^{\f 1{r'}}\\
&\leq |\theta_3-\theta_2|^{-\f d{r'}} \|g\|_{L^r}\|f_1\|_{L^{r'}}\|f_2\|_{L^\infty}\|f_3\|_{L^1},
\end{align*}
where in the last inequality we use Young's inequality,  in view of the fact that    $p_1=r'$ and $p_3=1$ for the  point $V_1^{(3)}$,   which implies $1/p_1+1/p_3=1+1/r'$. Note that $r$ could be equal to infinity in this case.
Boundedness for the  remaining
vertices follows by symmetry.

The case $m\ge 4$ is proved via the same method. Let $E^{(m)}$ be the intersection of the boundary     of the cube $[0,1]^m$  in $\mathbb R^m$ with the hyperplane $x_1+\cdots + x_m=1+1/r'$. Let us first consider the case where $r'>1$. The set $E^{(m)}$ is a closed convex polygon with $m(m-1) $ vertices which we call $V_i^{(m)}$. These are the   points
 in $\mathbb R^m$ whose coordinates are all zero except for two of them which are    $1$  and $1/r'$.
We claim that $E^{(m)}$ is the closed convex hull of the vertices $V_i^{(m)}$. In fact, we   show this by induction, as follows:
The assertion is true when  $m=3$. If this assertion is valid when $m=n-1$ for some $n\ge 4$, then
the faces of $E^{(n)}$ are the sets
$$
\{(1/q_1, \dots , 1/q_n):\,\, \sum_{i=1}^n q_i^{-1} = 1+ 1/r' \, , \quad \textup{$1\le q_i\le \infty$,  and   $q_i=\infty$ for only one $i$} \}.
$$
 Consider for instance the face
$$
F_1^{(n)}=\{(1/q_1, \dots , 1/q_n):\,\, \sum_{i=1}^n q_i^{-1} = 1+ 1/r' \, , \quad \textup{$1\le q_i\le \infty$ and   $q_1=\infty$} \}.
$$
 Then $F_1^{(n)}$ can be identified with the set
 $$
 \{(1/q_2, \dots , 1/q_n):\,\, \sum_{i=2}^n q_i^{-1} = 1+ 1/r' \, ,  \quad \textup{$1\le q_i\le \infty$} \}
 $$
which is equal to the intersection of the cube $[0,1]^{m-1}$ with the hyperplane $\sum_{i=2}^n x_i  = 1+ 1/r'$ in
$\mathbb R^{m-1}$
and  corresponds to the case $m=n-1$. Applying the    induction
hypothesis, the vertices of $F_1^{(n)}$ are the points $(1/q_1, \dots , 1/q_n)$  where $1/q_1=0$ and the
$1/q_i$ with $i\ge 2$ are   zero except for two of them which are either $1$ or $1/r'$. Thus, since $1/q_1=0$,  the vertices of
$F_1^{(n)}$ are the points $(1/q_1, \dots , 1/q_n)$  where $1/q_1=0$ and
$1/q_i$ with $i\ge 2$ are   zero except for two of them which are either
$1$ or $1/r'$. The same holds for the remaining faces and thus
the assertion that the points $V_i^{(m)}$ are the vertices of  $E^{(m)}$ holds when $m=n$. The same conclusion holds when
$r'=1$, but in this case the  intersection of the hyperplane $x_1+\cdots + x_m=2$ with the cube $[0,1]^m$  is a closed convex set with  $ m(m-1)/2$ vertices $W_i^{(m)}$ instead of $m(m-1)$ vertices.  Observe that the set of
$W_i^{(m)}$ is the set of all vectors with in $\mathbb R^m$ whose coordinates are zero except for two of them which are equal to $1$.

\begin{comment}
since the set $E$ is a closed convex polygon, there exists a line $l$ such that two vertices of set $E$ lie on $l$, i.e., one edge of $E$ is on $l$, whereas other vertices of set $E$ lie on the same side of $l$. Next, we pick one of two directions of $l$ as its positive direction. Assume all the vertices of set $E$ that not on $l$ are on the left side of $l$ according to this positive direction. Find the latter vertex, which we assume is $V_1$, of $E$ that lies on $l$ according to the positive direction. Then we rotate $l$ around $V_1$ anti-clockwise until $l$ meets other vertex of $E$. Assume this vertex is $V_2$, and we again rotate $l$ around $V_2$ anti-clockwise. Continuing these steps repeatedly until $l$ meets $V_1$ when rotating $l$ around $V_{m(m-1)}$ anti-clockwise. Notice that the step could be ended since $E$ has finite vertices. Thus, convex polygon formed by points $V_1,\dots, V_{m(m-1)}$ is the closed convex hull of vertices $V_1,\dots, V_{m(m-1)}$, while this convex polygon is the set $E$. Therefore, we prove our claim.
\end{comment}

 Notice that at each vertex $V_i^{(m)}$ (or $W_i^{(m)}$)
 boundedness holds and is proved
exactly in the same way as when $m=3$ with the only exception being that  the $L^\infty$ norm  of one function  is replaced by
the product of $L^\infty$ norms of the remaining $m-2$ functions. Then multilinear complex interpolation yields
boundedness in the closed convex hull of these vertices, i.e., yields
the required conclusion when $p=\frac{r'}{r'+1}$ (equivalently $q=1$).

Finally, the case where $\frac{r'}{1+r'}<  p< 1$ (equivalently $1< q< r$)  is obtained by
interpolation between the endpoint cases where
$p=\frac{r'}{r'+1}$ (equivalently $q=1$) and $p=1$ (equivalently $q=r$).   Notice that in the case where $p=1$ ($q=r$) the
intersection of the boundary of the unit cube $[0,1]^m$ with the hyperplane $x_1+\cdots + x_m=1$ is the
closed convex hull of the $m$ points $U_j^{(m)}=(0,\dots , 1, \dots , 0)$ with $1$ in the $j$th coordinate and zero in every other coordinate.
Let $Z_k^{(m)}$ be the vertices of the closed convex set formed by the intersection of the boundary of the unit cube $[0,1]^m$ with the hyperplane $x_1+\cdots + x_m=1/q+1/r'$ in $\mathbb R^m$ when $1<q<r$. Then there is a $\theta \in (0,1)$
[in fact $\theta = r'(\f 1q -\f 1r)$] such that
the set of all $Z_k^{(m)}$ is contained in the set
$$
\Big\{(1-\theta) U_j^{(m)} + \theta V_i^{(m)}{:} \qquad  j\in \{1,\dots ,m\}, \ \ i \in \{1,\dots , m(m-1)\}\Big\},
$$
where $V_i^{(m)}$ is replaced by $W_i^{(m)}$ for $i\in \{1,\dots , m(m-1)/2\}$,  if $r=\infty$.  Thus interpolation is possible in this case.

\section{Proof of Theorem \ref{thm1} }

\begin{proof} Using the inverse Fourier transform we write
\begin{align*}
&T_{\ell}(f_1,\dots,f_m)(x)=\int_{\R^d}f_1(x-\theta_1y)\cdots f_m(x-\theta_my)K_{\ell}(y)dy\\
&\ \ \ \ \ \ \ \ \ \ \ \ \ \ \ \ \ \ \ \ \ \ =\int_{\R^d}f_1(x-\theta_1y)\cdots f_m(x-\theta_my)e^{2\pi i\ell\cdot y}K(y)dy.
\end{align*}

We claim that the following inequality holds for $T$:
\begin{equation}\aligned \label{2.1}
T(f_1,\dots,f_m)(x)&=\Big(\sum\limits_{\ell\in \Z^d}|T_{\ell}(f_1,\dots,f_m)(x)|^2\Big)^{\frac{1}{2}}\\
&\leq \sum\limits_{u\in \Z^d}\bigg(\int_{Q_u}|f_1(x-\theta_1y)\cdots f_m(x-\theta_my)K(y)|^2dy\bigg)^{\frac{1}{2}}.
\endaligned\end{equation}
To verify the validity of this claim, we   take a sequence $a=\{a_{\ell}\}_{\ell\in\Z^d}\in \ell_2(\Z^d)$ with the property $\|a\|_{\ell_2(\Z^d)}=1$. Then by duality,  we need to prove   that
$$
\Big|\sum\limits_{\ell\in\Z^d}a_{\ell}T_{\ell}(f_1,\dots,f_m)(x)\Big|\leq \sum\limits_{u\in \Z^d}\bigg(\int_{Q_u}|f_1(x-\theta_1y)\cdots f_m(x-\theta_my)K(y)|^2dy\bigg)^{\frac{1}{2}}.
$$
Note that
\begin{align*}
&\sum\limits_{\ell\in\Z^d}a_{\ell}T_{\ell}(f_1,\dots,f_m)(x)=\int_{\R^d}f_1(x-\theta_1y)\cdots f_m(x-\theta_my)K(y)\sum\limits_{\ell\in\Z^d}a_{\ell}
e^{2\pi i\ell\cdot y}dy\\
&\ \ \ \ \ \ \ \ \ \ \ \ \ \ \ \ \ \ \ \ \ \ \ \ \ \ \ \ \ \ \ =\int_{\R^d}f_1(x-\theta_1y)\cdots f_m(x-\theta_my)K(y)\hat{a}(y)dy{ ,}
\end{align*}
where $\hat{a}$ is the Fourier transform of sequence $a$. It is easy to see that $\hat{a}$ is a periodic function with $\|\hat{a}\|_{L^2([0,1]^d)}=1$. Hence,
\begin{align*}
&\Big|\sum\limits_{\ell\in\Z^d}a_{\ell}T_{\ell}(f_1,\dots,f_m)(x)\Big|\\
&\leq \int_{\R^d}|f_1(x-\theta_1y)\cdots f_m(x-\theta_my)K(y)\hat{a}(y)|dy\\
& =\sum\limits_{u\in\Z^d}\int_{Q_u}|f_1(x-\theta_1y)\cdots f_m(x-\theta_my)K(y)\hat{a}(y)|dy\\
&\leq \sum\limits_{u\in\Z^d}\Big(\int_{Q_u}|f_1(x-\theta_1y)\cdots f_m(x-\theta_my)K(y)|^2dy\Big)^{1/2}\bigg(\int_{Q_u}|\hat{a}(y)|^2dy\bigg)^{1/2}\\
& =\sum\limits_{u\in\Z^d}\bigg(\int_{Q_u}|f_1(x-\theta_1y)\cdots f_m(x-\theta_my)K(y)|^2dy\bigg)^{1/2}.
\end{align*}
Therefore, we proved the claimed inequality (\ref{2.1}).

To show inequality (\ref{1.2}), we consider two cases:

\noindent\textbf{Case 1: $  2\le p<\infty$}.
Making use of Minkowski's inequality  we obtain
\begin{align*}
&\|T(f_1,\dots,f_m)\|_{L^p(\R^d)}\\
&\leq\sum\limits_{u\in\Z^d}\Big(\int_{\R^d}\Big(\int_{Q_u}|f_1(x-\theta_1y)\cdots f_m(x-\theta_my)K(y)|^2dy\Big)^{p/2}dx\Big)^{1/p}\\
&\leq\sum\limits_{u\in\Z^d}\Big(\int_{Q_u}\Big(\int_{\R^d}|f_1(x-\theta_1y)\cdots f_m(x-\theta_my)K(y)|^pdx\Big)^{2/p}dy\Big)^{1/2}\\
&=\sum\limits_{u\in\Z^d}\Big(\int_{Q_u}|K(y)|^2\Big(\int_{\R^d}|f_1(x-\theta_1y)\cdots f_m(x-\theta_my)|^pdx\Big)^{2/p}dy\Big)^{1/2}\\
&\leq\sum\limits_{u\in\Z^d}\Big\{\int_{Q_u}|K(y)|^2\Big[\Big(\int_{\R^d}|f_1|^{p\cdot\frac{p_1}{p}}dx\Big)^{\frac{p}{p_1}}
\cdots\Big(\int_{\R^d}|f_m|^{p\cdot\frac{p_m}{p}}dx\Big)^{\frac{p}{p_m}}\Big]^{\frac{2}{p}}dy\Big\}^{\frac{1}{2}}\\
&=B_2\prod\limits_{i=1}^m\|f_i\|_{L^{p_i}(\R^d)}.
\end{align*}

\noindent\textbf{Case 2: $1\le  p<2$}. By Minkowski's inequality, we obtain
\begin{align*}
&\|T(f_1,\dots,f_m)\|_{L^p(\R^d)}\\
&\leq \Big(\int_{\R^d}\Big|\sum\limits_{u\in\Z^d}\Big(\int_{Q_u}|f_1(x-\theta_1y)\cdots f_m(x-\theta_my)K(y)|^2dy\Big)^{1/2}\Big|^pdx\Big)^{1/p}\\
&\leq \sum\limits_{u\in\Z^d}\Big(\int_{\R^d}\Big(\int_{Q_u}|f_1(x-\theta_1y)\cdots f_m(x-\theta_my)K(y)|^2dy\Big)^{p/2}dx\Big)^{1/p}\\
&\leq\sum\limits_{u\in\Z^d}\Big\{\sum\limits_{n\in\Z^d}\Big(\int_{P_n}\int_{Q_u}|f_1(x-\theta_1y)\cdots f_m(x-\theta_my)K(y)|^2dydx\Big)^{p/2}\Big\}^{1/p}{ ,}
\end{align*}
where $P_n= n+[0, 1)^d$ and we used H\"older's inequality for exponents $2/p$ and $(2/p)'$ over the space $(P_n, dx)$ in the last inequality.

Let
$$
A_{n,u}=\int_{P_n}\int_{Q_u}|f_1(x-\theta_1y)\cdots f_m(x-\theta_my)K(y)|^2dydx
$$
and set   $f_{i}^{n,u} =f_i \chi_{P_n-\theta_i Q_u} $, $i=1,2,\dots,m$, $k_u =K \chi_{Q_u} $ and $\vec{f^{n,u}}=  (f_{i}^{n,u})_{1\le i\le m}$.
Using Theorem \ref{lem11} (part (b)) with $q=1$, $r=p'/2$, and $p_i/2$ in place of $p_i$, we obtain
\begin{align*}
A_{n,u}&\leq\int_{\R^d}\int_{\R^d}|f_{1}^{n,u}(x-\theta_1y)\cdots f_{m}^{n,u}(x-\theta_my)k_u(y)|^2dydx\\
&= \int_{\R^d}|\vec{f^{n,u}}|^2\otimes|k_u|^2(x)dx \\
%& \leq \||\vec{f_n}|^2\otimes|k_u|^2\|_{L^{q/2}(\R^d)}\\
&\leq C\prod\limits_{i=1}^m\| |f_{i}^{n,u}|^2\|_{L^{p_i/2}(\R^d)}\| |k_u|^2\|_{L^{p'/2  }(\R^d)}\\
&=C\prod\limits_{i=1}^m\|f_{i}^{n,u}\|_{L^{p_i}(\R^d)}^2\|  k_u \|_{L^{p'   }(\R^d)}^2.
\end{align*}
Therefore, it follows from H\"older's inequality that
\begin{align*}
&\|T(f_1,\dots,f_m)\|_{L^p(\R^d)}\\
&\leq C\sum\limits_{u\in\Z^d}\Big(\sum\limits_{n\in\Z^d}\Big(\prod\limits_{i=1}^m
\|f_{i}^{n,u}\|_{L^{p_i}(\R^d)}^2\|k_u\|^2_{L^{p'  }(\R^d)}\Big)^{p/2}\Big)^{1/p}\\
&= C\sum\limits_{u\in\Z^d}\|k_u\|_{L^{p'   }(\R^d)}\Big(\sum\limits_{n\in\Z^d}
\Big(\int_{\R^d}|f_{1}^{n,u}|^{p_1}dx\Big)^{\frac p{p_1}}\cdots\Big(\int_{\R^d}|f_{m}^{n,u}|^{p_m}dx
\Big)^{\frac p{p_m}}\Big)^{\frac 1p}\\
&\leq C\sum\limits_{u\in\Z^d}\|K \chi_{Q_u}\|_{L^{p'  }(\R^d)} \Big\{\prod_{j=1}^m \Big(\sum\limits_{n\in\Z^d}
\Big(\int_{P_n-\theta_j Q_u}|f_j|^{p_j}dx\Big)^{\frac{p}{p_j}\frac{p_j}{p}}\Big)^{\frac{p}{p_j}} \Big\}^{\frac{1}{p}}\\
&\leq C\, B_p   \prod\limits_{i=1}^m\|f_i\|_{L^{p_i}(\R^d)} .
\end{align*}

Thus, we complete the proof of Theorem \ref{thm1}.

\par\bigskip
\end{proof}

%%%%%%%%%%%%%%%%%% References %%%%%%%%%%%%%%%%%%%%%%%%%%%%%%%%%%%%%%

\end{document}